\documentclass{article}
\usepackage[english]{babel}
\usepackage[utf8]{inputenc}

\usepackage{bm}
\usepackage{amsmath}
\usepackage{graphicx}
\usepackage{color}
\usepackage{multirow}
\usepackage{amsfonts}
\newcommand{\grad}{\mathop{\rm grad}\nolimits}
\renewcommand{\div}{\mathop{\rm div}\nolimits}

\graphicspath{ {./} } 

\date{}

\author{
Aleksei Tyrylgin
\thanks{Multiscale model reduction laboratory, North-Eastern Federal University, Yakutsk, Republic of Sakha (Yakutia), Russia, 677980.}
\and
Maria Vasilyeva \thanks{Institute for Scientific Computation, Texas A\&M University, College Station, TX 77843-3368 \& Department of Computational Technologies, North-Eastern Federal University, Yakutsk, Republic of Sakha (Yakutia), Russia, 677980. Email: {\tt vasilyevadotmdotv@gmail.com}.}
\and
Denis Spiridonov
\thanks{Multiscale model reduction laboratory, North-Eastern Federal University, Yakutsk, Republic of Sakha (Yakutia), Russia, 677980.}
\and
Eric T. Chung \thanks{Department of Mathematics,
The Chinese University of Hong Kong (CUHK), Hong Kong SAR. Email: {\tt tschung@math.cuhk.edu.hk}.}
}

\begin{document}
\title{Generalized Multiscale Finite Element Method for the poroelasticity problem in multicontinuum media}

\maketitle

\begin{abstract}
In this paper, we consider a poroelasticity problem in heterogeneous multicontinuum media that is widely used in simulations of the unconventional hydrocarbon reservoirs and geothermal fields.
Mathematical model contains a coupled system of equations for pressures in each continuum and effective equation for displacement with volume force sources that are proportional to the sum of the pressure gradients for each continuum.
To illustrate the idea of our approach, we consider a dual continuum background model with discrete fracture networks that can be generalized to a multicontinuum model for poroelasticity problem in complex heterogeneous  media.
We present a fine grid approximation based on the finite element method and Discrete Fracture Model (DFM) approach for two and three-dimensional formulations.
The coarse grid approximation is constructed using the Generalized Multiscale Finite Element Method (GMsFEM), where we solve local spectral problems for construction of the multiscale basis functions for displacement and pressures in multicontinuum media.
We present numerical results for the two and three dimensional model problems in heterogeneous fractured porous media. We investigate relative errors between reference fine grid solution and presented coarse grid approximation using GMsFEM with different numbers of multiscale basis functions.
Our results indicate that the proposed method is able to give accurate solutions with few degrees of freedoms.

\end{abstract}

\section{Introduction}
Mathematical models of the poroelasticity problems in multicontinuum media for heterogeneous fractured media are used for simulation of the unconventional hydrocarbon reservoirs, geothermal fields, underground disposal of radioactive waste in subsurface collectors, etc \cite{wu2011multiple, li2018multiscale, dpgmsfem2017, akkutlu2012multiscale, wu2007triple}.
Mathematical model is described by a system of equations for pressures in each continuum and equation for displacements.
The most important feature of the mathematical model is that the equations are coupled. We can highlight two main coupling types: (1) multicontinuum pressures coupling via mass transfer term \cite{barenblatt1960basic, warren1963behavior, arbogast1990derivation}, and (2) pressure and displacement coupling via term that describes the compressibility of the medium and volume force which is proportional to the pressure gradient \cite{kolesov2014splitting, kim2010sequential}.

Mathematical models for flow in the multicontinuum media is used for describing of the complex flow processes in multiscale fractured heterogeneous porous media \cite{martin2005modeling, formaggia2014reduced, d2012mixed, Quarteroni2008coupling}.
The flow in fractures has a significant impact on filtration processes and requires careful consideration \cite{lee2001hierarchical, li2006efficient}. Moreover, since the fractures are characterized by high permeability and their thickness is significantly smaller than the size of the simulated field, this leads to the need to build special mathematical models of multicontinuum, where independent variables are distinguished to describe the flow in a porous medium and in the network of fractures taking into account the special flow function \cite{pororelgmsfem2018, li2008efficient, karimi2003efficient}.
At the same time, the size of fractures should be separated, since they can exist on different scales and can differ in the nature of their occurrence. In the case of naturally fractured porous media, the fracture system is basically connected and dual porosity models are traditionally used \cite{barenblatt1960basic, warren1963behavior}.
The interaction of the continua is described by specifying the flow functions between continua (mass transfer) \cite{tene2015algebraic, li2008efficient, dpgmsfem2017}.

For numerical simulation of the problems in fractured and heterogeneous porous media, we should use a sufficiently fine grid that resolve all small scale features in the level of mesh construction. Therefore, the discrete formulation of such problems leads to the large system of equations that is computationally expensive. To reduce the size of the discrete system, various multiscale methods have been developed \cite{egh12gmsfem, gmsfem2016adaptive, eh09msfem, jenny2005adaptive, lunati2006multiscale, hajibeygi2008iterativefvm}.
In \cite{hkj12, ctene2016algebraic}, the multiscale finite volume method was presented for solution of the flow problems in fractured porous media. Multiscale finite volume method for solution of the poroelasticity problem is presented in \cite{castelletto2018multiscale}. Generalized multiscale finite element (GMsFEM) for solution of the flow problems in fractured porous media is considered in \cite{efendiev2015hierarchical, akkutlu2015multiscale}. In \cite{brown2016generalized1, brown2016generalized2}, we considered construction of the coarse grid problem for poroelasticity problem in heterogeneous media. Extension of the GMsFEM method for solution of the poroelasticity problem in fractured media with discontinuous Galerkin method for displacements and continuous Galerkin method for pressure is presented in our previous work \cite{pororelgmsfem2018}. Recently a new method was presented for solution of the flow and poroelasticity problems in fractured and heterogeneous porous media \cite{nlmc2018, chung2017cem, poroelnlmc2018}. Generalization of the GMsFEM and NLMC approach for solution of the flow problems in multicontinuum media is considered in \cite{dpgmsfem2017, vasilyeva2019nlmcdp}.

In this paper, we consider the Generalized Multiscale Finite Element method for solution of the poroelasticity problems in multicontinuum heterogeneous media.
We construct a multiscale basis functions for pressures and displacement. Construction of the basis functions for flow problem in multicontinuum media is based on the solution of the coupled system of equations in each local domains. Method automatically identify each continuum flow features via solution of the local spectral problems. For displacement, we use similar approach where main features are captured by the local spectral problem. We numerically investigate presented method for two and three-dimensional model problems in fractured and heterogeneous porous media.

The work is organized as follows. In Section 2, we present the mathematical model of the poroelasticity problem in multicontinuum medium. Then in Section 3, a fine grid approximation is constructed using the finite element method. In Section 4, we present a coarse grid approximation using Generalized Multiscale Finite Element method, where we describe construction of the multiscale basis functions and coarse grid system construction. Numerical results for two and three-dimensional model poroelasticity problems are presented in Section 5. Finally, we present Conclusions.

\section{Mathematical model}

Let $\Omega \subset R^d$ is the computational domain for background medium with dual continuum approach, where $d = 2$ for two-dimensional problems and $d = 3$ for three-dimensional problems.
For example, the first continuum can describe a flow in the matrix of the porous media, and the second continuum belongs to the network of small highly connected fracture network.
Furthermore, we define $\gamma \in R^{d-1}$ as a computational domain for low dimensional fracture networks model that describes the flow in the large-scale fractures.
For the fluid flow in the poroelastic medium, we have the mass balance equation and the Darcy's law
\begin{equation}
\begin{split}
\alpha_1 \frac{\partial \div u}{\partial t}
+ \frac{1}{M_1} \frac{\partial p_1}{\partial t}
+ \div q_1 + q_{12} + q_{1f} = f_1, \ \ x \in \Omega, \\
q_1 = -k_1 \grad p_1, \ \ x \in \Omega, \\
\alpha_2 \frac{\partial \div u}{\partial t}
+ \frac{1}{M_2} \frac{\partial p_2}{\partial t}
+ \div q_2 - q_{21} + q_{2f} = f_2, \ \ x \in \Omega, \\
q_2 = -k_2 \grad p_2, \ \ x \in \Omega, \\
\alpha_f \frac{\partial \div u}{\partial t}
+ \frac{1}{M_f} \frac{\partial p_f}{\partial t}
+ \div q_f - q_{f1} - q_{f2} = f_f, \ \ x \in \gamma, \\
q_f = -k_f \grad p_f, \ \ x \in \gamma,
\end{split}
\end{equation}
where $u$ is the displacement, $q_i$ is the velocity, $p_i$ is the pressure, $\kappa_i$ is the permeability ($k_i = \kappa_i / \mu$, $\mu$ is the fluid viscosity), $f_i$ refer to source and sink terms, $\alpha_i$ is Biot coefficients, $M_i$ is the Biot modulus and $i$ is continuum index ($i = 1,2,f)$.

Here $q_{12}$, $q_{1f}$ and $L_{2f}$ are the transfer term between first and second continua; fist continuum and fractures; and second continuum and fractures, respectively.
We have
\[
q_{12} = q_{21} = r_{12}(p_1 - p_2), \quad
r_{if} = \eta_i r_{if} (p_i - p_f), \quad
r_{fi} = \eta_f r_{if} (p_i - p_f),
\]
and $\int_{\Omega} q_{12} dx = \int_{\Omega} q_{21} dx $, $\int_{\Omega} q_{if} dx = \int_{\gamma} q_{fi} ds$, $\eta_i$ is the geometric factor and $r_{ij}$ is the mass transfer term that proportional to the continuum permeabilities.

For the mechanics of the poroelastic multicontinuum media, we use an effective equation for displacement with volume force sources that proportional to the sum of the pressure gradients for each continuum.
\begin{equation}
\begin{split}
- \div \sigma_T (u,p_1, p_2, p_f) = 0, \quad  x \in \Omega, \\
\sigma_T(u, p_1, p_2, p_f) = \sigma(u) - \sum_{i= 1,2,f} \alpha_i p_i \mathcal{I}, \quad  x \in \Omega,
\end{split}
\end{equation}
where $\sigma_T$ is the total stress tensor, $\sigma$ is the stress tensor  \cite{pororelgmsfem2018, poroelnlmc2018, kim2010sequential}. 
In the case of a linear elastic stress-strain constitutive relation, we have
\[
\sigma(u) = 2 \mu \varepsilon(u) + \lambda \div u \, \mathcal{I}, \ \  \varepsilon(u) = \frac{1}{2}(\grad u + \grad u^T),
\]
where $\varepsilon$ is the strain tensor, $\lambda$ and $\mu$ are the Lame’s coefficients.

Then, we have the following coupled system of equations 
\begin{equation}
\label{eq:main}
\begin{split}
 \alpha_1 \frac{\partial \div  u}{\partial t} + \frac{1}{M_1} \frac{\partial p_1}{\partial t} - \div \cdot (k_1 \grad p_1)  + q_{12} +  q_{1f}  = f_1, 
 \quad x \in  \Omega, \\
 \alpha_2 \frac{\partial \div  u}{\partial t} + \frac{1}{M_2} \frac{\partial p_2}{\partial t} - \div \cdot (k_2 \grad p_2)  - q_{21} + q_{2f}  = f_2, 
 \quad x \in \Omega, \\ 
 \alpha_f \frac{\partial \div  u}{\partial t} +  \frac{1}{M_f} \frac{\partial p_f}{\partial t} - \div \cdot (k_f \grad p_f)  - q_{f1} - q_{f2}  = f_f , 
 \quad x \in \gamma, \\
 - \div  \sigma (u) + \alpha_1 \grad p_1  + \alpha_2  \grad p_2 + \alpha_f \grad p_f  = 0 , 
 \quad x \in \Omega. 
\end{split}
\end{equation}

We consider a system of equations \eqref{eq:main} with the following initial conditions 
\begin{equation}
\label{eq:main-ic}
p_1 =  p _2 = p_f = p^0, \quad u = 0, \quad  x \in \Omega,
\end{equation}
and boundary conditions
\begin{equation}
\label{eq:main-bc}
\begin{split}
p_f = g, \quad x \in \Gamma_{\gamma} \quad  
-k_f \frac{\partial p_f}{\partial n} = 0, \quad x \in \partial \gamma / \Gamma_{\gamma}, \\
-k_1 \frac{\partial p_1}{\partial n} = 0, \quad  x \in \Omega,  \quad
-k_2 \frac{\partial p_2}{\partial n} = 0, \quad x \in \Omega, \\
u = 0, \quad x \in \Gamma_D \quad 
\sigma = 0, \quad x \in \Gamma_N,
\end{split}
\end{equation}
where boundaries 
$\Gamma_D \cup \Gamma_N = \partial \Omega$ and $\Gamma_{\gamma} = \partial \Omega \cap \partial \gamma$. 

We can generalize presented model as poroelasticity model for multicontinuum media
\begin{equation}
\label{eq:main-mc}
\begin{split}
 \alpha_i \frac{\partial \div  u}{\partial t} 
 + \frac{1}{M_i} \frac{\partial p_i}{\partial t} 
 - \div \cdot (k_i \grad p_i)  + \sum_{j \neq i}  r_{ij} (p_i - p_j)  = f_i, 
 \quad x \in  \Omega, \\
 - \div  \sigma (u) + \sum_{j} \alpha_j \grad p_j  = 0 , 
 \quad x \in \Omega. 
\end{split}
\end{equation}
where $i = 1,...,M$ and $M$ is the number of continua.

\section{Fine grid finite element approximation}

For the approximation of the system of equations \eqref{eq:main} with boundary conditions \eqref{eq:main-bc}, we use a finite element method. 
To do so, we define the following functional spaces
\[
V = \lbrace v \in H^1(\Omega):  v = 0 \text{ on } \Gamma_D \rbrace,  \quad 
W_1 = W_2 =  H^1(\Omega), 
\]\[
W_f = \lbrace 
w \in H^1(\gamma): w = g \text{ on } \Gamma_{\gamma} 
\rbrace, \quad 
\hat{W}_f = \lbrace 
w \in H^1(\gamma): w = 0 \text{ on } \Gamma_{\gamma} 
\rbrace.
\]

The 
variational formulation of  
the poroelasticity problem in multicontinuum media can be written as follows: 
find $(p_1, p_2, p_f, u) \in W_1 \times W_2 \times W_f \times V$ such that
\begin{equation}
\label{eq:fine}
\begin{split}
& d_i \left( \frac{\partial u}{\partial t}, w_i \right) 
+ c_i \left( \frac{\partial p_i}{\partial t}, w_i \right) 
+ b_i(p_i, w_i) 
+ \sum_j q_{ij}(p_i - p_j, w_i) 
=  l(w_i), 
\quad \forall w_i \in W_i, 
\\
& a(u, v) + \sum_j g_i(p_i, v) = 0, 
\quad \forall v \in V, 
\end{split}
\end{equation}
where bilinear and linear forms are following
\[
b_i(p_i, w_i) = \int_{\Omega_i} k_i \nabla p_i \cdot \nabla w_i dx, \quad
l(w_i) = \int_{\Omega_i} f_i w_i dx,
\]\[
c_i(p_i, w_i) =  \int_{\Omega_i} \frac{1}{M_i} p_i w_i \, dx, \quad 
d_i(u, w_i) =  \int_{\Omega_i} \alpha_i \div u \, w_i \, dx,
\]\[
q_{ij}(p_i - p_j, w_i) =  \int_{\Omega_i} r_{ij}(p_i - p_j) \, w_i \, dx,
\]\[
a(u,v) = \int_{\Omega} \sigma(u) \cdot \varepsilon(v) \, dx, \quad 
g_i(p_i, v) =  \int_{\Omega_i} \alpha_i \grad p_i v \, dx,
\]
and $\Omega_1 = \Omega_2 = \Omega$, $\Omega_f = \gamma$, $i = 1,2,f$.

Let $\mathcal{T}^h$ denote a finite element partition of the domain $\Omega$. 
For the fracture continuum, we use a discrete fracture model and use an unstructured fine grid $\mathcal{T}^h$ that explicitly resolve fracture geometry.
We assume that $\cup_j \gamma_j$ is the subset of faces for $\mathcal{T}^h$ that represent fractures, where $j = 1,...,N_{frac}$ and $N_{frac}$ is the number of discrete fractures. For approximation by time, we use an implicit finite difference scheme with time step $\tau$. Therefore, we have following discrete system in matrix form on the fine grid for the triple-continuum media and $y^h = (p_1^h, p_2^h, p_f^h, u^h)^T$
\begin{equation}
C \frac{y^h - \check{y}^h}{\tau} + A y^h = F, 
\end{equation}
where  
\[
C = 
\begin{pmatrix}
  C_1 & 0 & 0 & D_1 \\
  0 & C_2 & 0 & D_2 \\
  0 & 0 & C_f & D_f \\
  0 & 0 & 0 & 0
\end{pmatrix}, \quad 
A = 
\begin{pmatrix}
A_1+Q_{12} + Q_{1f} & -Q_{12} & -Q_{1f} & 0 \\
-Q_{12} & A_2+Q_{12} + Q_{2f} & -Q_{2f} & 0 \\
-Q_{1f} & -Q_{2f} & A_f+Q_{1f} + Q_{2f} & 0 \\
D_1^T & D_2^T & D_f^T & A_u
\end{pmatrix}
\]
and $F = (F_1, F_2, F_f, 0)^T$, $\check{y}^h$ is the solution from the previous time step. 
Here 
\[ 
A_i = [a_{i,ln}], \quad 
a_{i, ln} = \int_{\Omega_i} k_i \grad \phi^i_l \cdot \grad \phi^i_n dx, \quad
A_u = [a_{u,ln}], \quad  
a_{u,ln} = \int_{\Omega} \sigma(\Phi_l) \cdot \varepsilon(\Phi_n) dx,
\]\[
Q_{ij} = [q_{ij, ln}], \quad 
q_{ij, ln} = \int_{\Omega_i} r_{ij} \phi^i_l \phi^j_n dx, \quad 
C_i=[c_{i,ln}], \quad 
c_{i,ln}=\int_{\Omega_i} \frac{1}{M_i} \phi^i_l \phi^i_n dx,  
\]\[ 
D_i = [d_{i,ln}], \quad 
d_{i,ln} = \int_{\Omega} \alpha_i  \div \Phi_l \phi^i_n  dx,  \quad
F_i = [f_{i,l}], \quad 
f_{i,l} = \int_{\Omega_i} f_i \phi^i_l dx,
\]
and $p_i = \sum_l p_{i,l}^h \phi^i_l$,  $u = \sum_l u_l^h \Phi_l$, where 
$\Phi_l$ is the linear basis functions for displacements, $\phi^1_l = \phi^2_l$ is the $d$ - dimensional linear basis functions for pressure, $\phi^f_l$ is the  $(d-1)$ - dimensional linear basis functions for pressure. 

In this paper, for simplification of the matrix construction, we use a modified DFM approach and consider the case when $\alpha_f = 0$, $\sigma_{2f} = 0$. We assume that $p^h_1= p^h_f$  and using superposition principle \cite{efendiev2015hierarchical, akkutlu2015multiscale}, we eliminate $p^h_f$ from equations and obtain following coupled system of equations for $y^h = (p_1^h, p_2^h, u^h)^T$
\[
\tilde{C} \frac{y^h - \check{y}^h}{\tau} + \tilde{A} y^h = \tilde{F},
\]
where
\[
\tilde{C} = 
\begin{pmatrix}
C_1+C_f & 0 & D_1 \\
0 & C_2  & D_2 \\
0 & 0  & 0
\end{pmatrix}, \quad 
\tilde{A} = 
 \begin{pmatrix}
A_1 + A_f+Q_{12} & -Q_{12} & 0 \\
 -Q_{12} & A_2+Q_{12} & 0\\
D^T_1 & D^T_2 & A_u
\end{pmatrix}, \quad 
\tilde{F} = 
 \begin{pmatrix}
F_1 + F_f \\
F_2 \\
0
\end{pmatrix}. 
\]
with matrices of the size $N_h = 4 N_v$,  $N_v$ is the number of vertices in $\mathcal{T}^h$.

\section{Coarse grid approximation using GMsFEM}

For construction of the coarse grid approximation of the poroelasticity problems in fractured and heterogeneous media, we use a Generalized Multiscale Finite Element Method (GMsFEM).
GMsFEM contains following steps:
\begin{enumerate}
\item coarse grid and local domains construction;
\item solution of the local problems with different boundary conditions to construct a snapshot space in each local domain;
\item multiscale basis functions construction via solution of the local spectral problems on the snapshot space;
\item generation of the projection matrix using local multiscale basis functions;
\item construction of the coarse grid system using projection matrix;
\item solution of the unsteady problem on the coarse grid and reconstruction of the fine grid solution.
\end{enumerate}
In this computational algorithm, first four steps are offline (preprocessing) steps for a given fracture geometry and heterogeneity. Fifth step is also offline for linear problems and time - independent right-hand side, but should be online step for nonlinear problems, where fine grid system is change on each nonlinear or/and time iteration. After that on the sixth (online) step, we can perform fast and accurate solution of the reduced order model on the coarse grid. 

In this work, we construct multiscale basis functions for displacements and pressures separately, but basis functions for multicontinuum pressure equations are constructed in the coupled way. We start with definition of the snapshot space, after that we define a local eigenvalue problems for pressures and displacements. Finally, we define projection matrix and present construction of the coarse grid poroelasticity system on the multiscale space.
Let $\mathcal{T}^H$ is the coarse grid partitioning of the domain
\[
\mathcal{T}^H = \bigcup_j K_j,
\]
where $K_j$ is the coarse grid cell. We will use a continuous Galerkin approximation on the coarse grid, and define local domain $\omega_l$ for multiscale basis functions as combination of the several coarse grid cells that share same coarse grid node ($l = 1,...,N_v^H$, $N_v^H$ is the number of coarse grid vertices).

\textbf{Multiscale basis functions for pressures in multicontinuum media.}
To construct a snapshot space, we solve following local problem in domain $\omega_l$:
find $\psi^{l, j} = (\psi^{l, j}_1, ..., \psi^{l, j}_M) \in W^h_1 \times ... \times W^h_M$ such that
\begin{equation}
b_i(\psi^{l, j}_i, w_i) + \sum_j q_{ij}(\psi^{l, j}_i - \psi^{l, j}_j, w_i)
= 0, \quad \forall w_i \in \hat{W}^h_i
\end{equation}
where
\[
W^h_i = \lbrace
w \in H^1(\omega_l): w = \delta^j_i \text{ on } \partial \omega_l
\rbrace, \quad
\hat{W}^h_i = \lbrace
w \in H^1(\omega_l): w = 0 \text{ on }\partial \omega_l
\rbrace,
\]
and $\delta^j_i$ is the piecewise constant function (delta function) for $j = 1,..,N_v^{\omega_l}$ ($,N_v^{\omega_l}$ is the number of nodes on the computation mesh for $\omega_l$), $i$ is the index of continuum ($i = 1,...,M$). Therefore, we solve $L_p^{\omega_l} = N_v^{\omega_l} \cdot M$ local problems.

We define snapshot space for pressures in multicontinuum media as follows.\begin{equation}
W_{snap}(\omega_l) = \text{span} \lbrace 
\psi^{l, j}, \, l = 1,...,N^H_v, \, j = 1,...,L_p^{\omega_l}
\rbrace.
\end{equation}
Next, we solve following local spectral problem on the snapshot space
\begin{equation}
\tilde{A}_p \tilde{\phi}^l = \lambda_p \tilde{S}_p \tilde{\phi}^l, 
\end{equation}
where  $\hat{\phi}^l = (R^p_{snap})^T \tilde{\phi}^l$ and 
\[
\tilde{A}_p = R^p_{snap} {A}_p (R^p_{snap})^T, \quad 
\tilde{S}_p = R^p_{snap} {S}_p (R^p_{snap})^T, \quad 
R^p_{snap} = (\psi^{l, 1},...,\psi^{l, L_p^{\omega_l}})^T
\]
Here for matrices in triple continuum case, we have
\[
S_p = 
\begin{pmatrix}
  S_1 & 0 & 0 \\
  0 & S_2 & 0 \\
  0 & 0 & S_f 
\end{pmatrix}, \quad 
A_p = 
\begin{pmatrix}
A_1+Q_{12} + Q_{1f} & -Q_{12} & -Q_{1f} \\
-Q_{12} & A_2+Q_{12} + Q_{2f} & -Q_{2f} \\
-Q_{1f} & -Q_{2f} & A_f+Q_{1f} + Q_{2f}
\end{pmatrix}
\]
where 
\[
A_i = [a_{i,mn}], \quad 
a_{i, mn} = \int_{\omega_l^i} k_i \grad \phi^i_m \cdot \grad \phi^i_n dx, \quad 
S_i=[s_{i,mn}], \quad 
s_{i,mn}=\int_{\omega_l^i} k_i \phi^i_m \phi^i_n dx. 
\]
We choose an eigenvector $\hat{\phi}_j$ ($j = 1,..,M^{l,p}$) corresponding to the first smallest $M^{l,p}$ eigenvalues and multiply to the linear partition of unity functions $\chi^l$ for obtaining conforming basis functions
\[
W_{ms} = \text{span} \lbrace 
\phi^{l,j}, \, l = 1,...,N^H_v, \, j = 1,...,M^{l,p}
\rbrace,
\]
where $\phi^{l,j} = \chi^l \hat{\phi}^{l,j}$.

\textbf{Multiscale basis functions for displacements.}
We construct the multiscale basis functions by solution following problem in local domain $\omega_l$:
find $\Psi^{l, j} \in V^h$ such that
\begin{equation}
a(\Psi^{l, j}, v) = 0,
\quad \forall v \in \hat{V}^h,
\end{equation}
where
\[
V^h = \lbrace
v \in H^1(\omega_l): v = \bar{\delta}^j_i \text{ on } \partial \omega_l
\rbrace, \quad
\hat{V}^h = \lbrace
v \in H^1(\omega_l): v = 0 \text{ on } \partial \omega_l
\rbrace.
\]
and $\bar{\delta}^j_i$ is the vector for each component for $d$ - dimensional problem ($d = 2,3$) i.e. $\bar{\delta}^j_i = ({\delta}^j_i, 0, 0)$ or $\bar{\delta}^j_i = (0, {\delta}^j_i, 0)$ or $\bar{\delta}^j_i = (0, 0, {\delta}^j_i)$ for $d=3$. We solve $L_u^{\omega_l} = d \cdot N_v^{\omega_l}$ local problems.

We define snapshot space for pressures in multicontinuum media as follows
\begin{equation}
V_{snap}(\omega_l) = \text{span} \lbrace 
\Psi^{l, j}, \, l = 1,...,N^H_v, \, j = 1,...,L_u^{\omega_l}
\rbrace.
\end{equation}
For construction of the multiscale basis, we solve following local spectral problem on the snapshot space
\begin{equation}
\tilde{A}_u \tilde{\Phi} = \lambda_u \tilde{S}_u \tilde{\Phi}, 
\end{equation}
where $\hat{\Phi}^l = (R^u_{snap})^T \tilde{\Phi}^l$, 
\[
\tilde{A}_u = R^u_{snap} {A}_u (R^u_{snap})^T, \quad 
\tilde{S}_u = R^u_{snap} {S}_u (R^u_{snap})^T, \quad 
R^u_{snap} = (\Psi^{l, 1},...,\Psi^{l, L_u^{\omega_l}})^T
\]
and 
\[
A_u = [a_{u,mn}], \quad  
a_{u,mn} = \int_{\omega_l} \sigma(\Phi_m) \cdot \varepsilon(\Phi_n) dx,
\quad 
S_u=[s_{i,mn}], \quad 
s_{i,mn}=\int_{\omega_l} (\lambda + 2\mu) \Phi^i_m \Phi^i_n dx. 
\]
We choose an eigenvector $\hat{\Phi}_j$ ($j = 1,..,M^{l,u}$) corresponding to the first smallest $M^{l,u}$ eigenvalues and multiply to the linear partition of unity functions for obtaining conforming basis functions
\[
V_{ms} = \text{span} \lbrace 
\Phi^{l,j}, \, l = 1,...,N^H_v, \, j = 1,...,M^{l,u}
\rbrace,
\]
where $\Phi^{l,j} = \chi^l \hat{\Phi}^{l,j}$.

\textbf{Coarse grid system}. 
Using constructed multiscale basis functions, we define projection matrix
\begin{equation}
R = \begin{pmatrix}
R_p & 0 \\
0 & R_u
\end{pmatrix}
\end{equation}
where 
\[
R_u = (\Phi^{1,1}, ... , \Phi^{1,M^{1,u}} , ..., \Phi^{N^H_v,1}, ...,\Phi^{N^H_v, M^{N^H_v,u}})^T,
\]\[
R_p = (\phi^{1,1}, ... , \phi^{1,M^{1,p}} , ..., \phi^{N^H_v,1}, ...,\phi^{N^H_v, M^{N^H_v,p}})^T.
\]
Finally, we obtain following reduced order model
\[
{C}^H \frac{y^H - \check{y}^H}{\tau} + {A}^H y^H = {F}^H,
\]
where $C^H = R C R^T$, $A^H = R A R^T$ and $F^H = R F$. 
with matrices of the size $N_f = 4 N_v$,  $N_v$ is the number of vertices in $\mathcal{T}^h$. 
Size of the system is $N_H = \sum_{l = 1}^{N^H_v} (M^{l,p} + M^{l,u}) $ or $N_H = (M^p + M^u) \cdot N^H_v$ for $M^p = M^{l,p}$ and $M^u = M^{l,u}$  ($\forall l = 1,...,N^H_v$), where $N^H_v$ is the number of vertices of coarse grid $\mathcal{T}^H$.

After obtaining of a coarse-scale solution, we reconstruct fine-scale solution 
\[
y^{ms} = R^T y^H.
\]
We note that, in general, multiscale basis functions for the displacements and pressures can be calculated by solution coupled poroelasticity problem similarly to the multicontinuum pressures.

\section{Numerical results}

\begin{figure}[h!]
\centering
\includegraphics[width=0.4\linewidth]{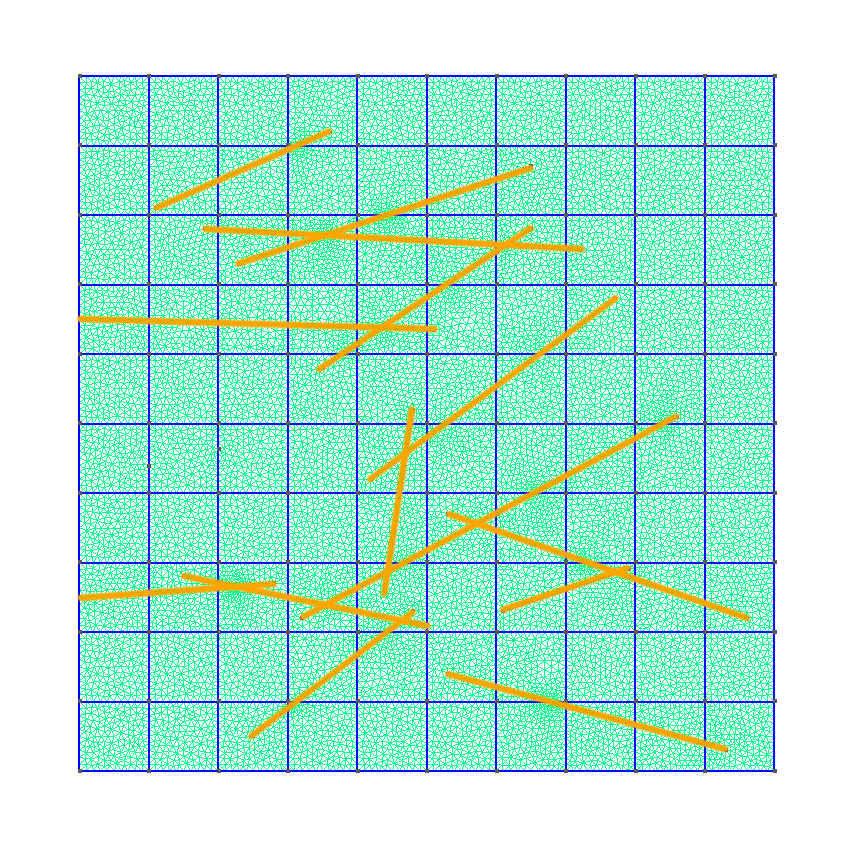}
\caption{Computation domain and grid for two - dimensional problem. Coarse grid (blue color), fine grid (green) and fractures (orange)}
\label{mesh2d}
\end{figure}

\begin{figure}[h!]
\centering
\includegraphics[width=0.75\linewidth]{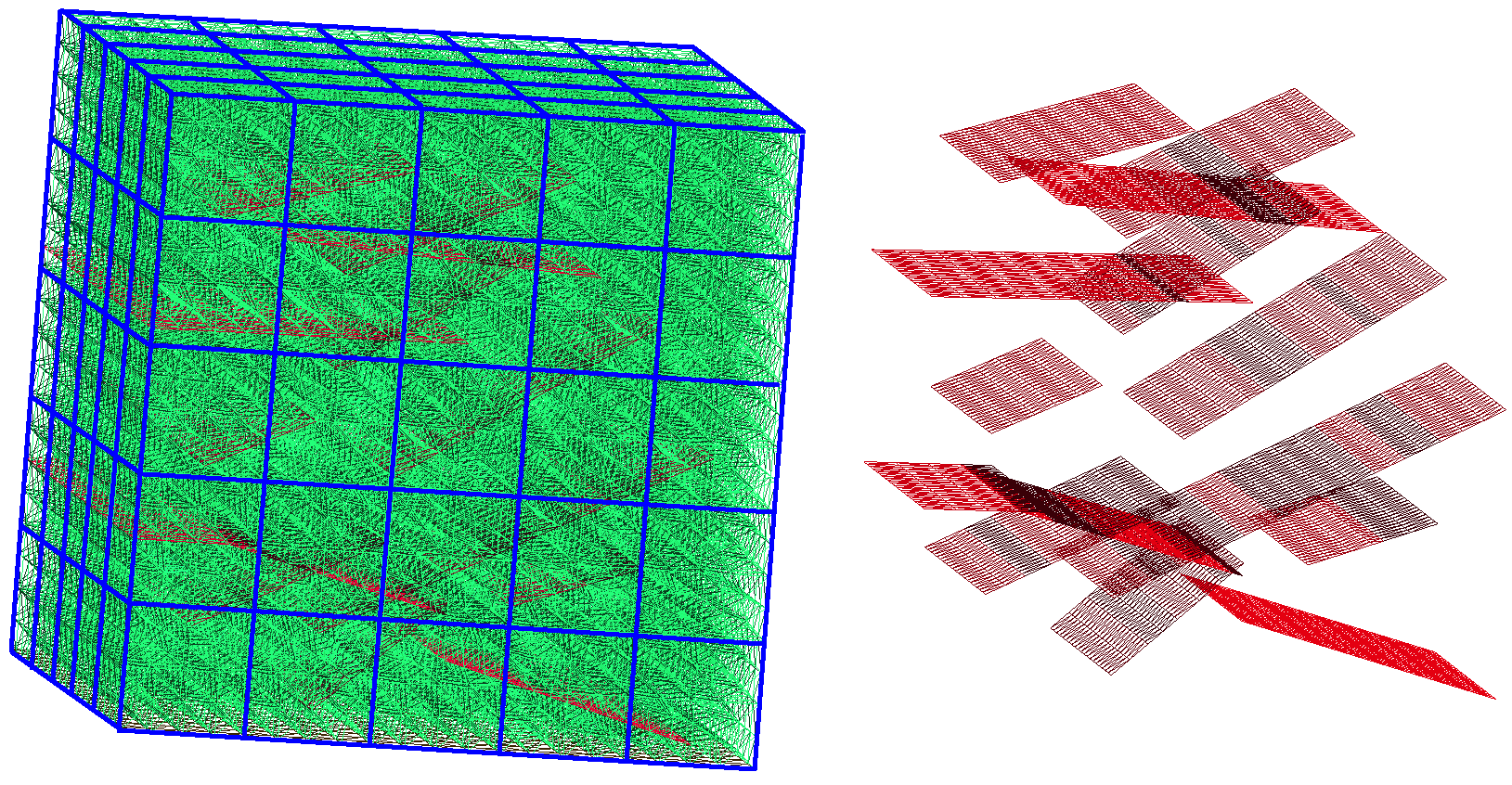}
\caption{Computation domain and grid for three - dimensional problem. Coarse grid (blue color), fine grid (green) and fractures (red)}
\label{mesh3d}
\end{figure}

In this section, we consider poroelasticity problem in fractured and heterogeneous media.
We consider two problems in $\Omega = [0, 10]^d$ ($d = 2,3$):
\begin{itemize}
\item \textit{Two - dimensional model problem}. Fine grid contains 14376 vertices and 28350 cells. Coarse grid contains 121 vertices and 100 cells (Figure \ref{mesh2d});
\item \textit{Three - dimensional model problem}. Fine grid contains 21609 vertices and 118500 cells. Coarse grid contains 216 vertices and 125 cells (Figure \ref{mesh3d}).
\end{itemize}

We set parameters of model problem as follows:
$\alpha_1 = 1$, $\alpha_2 = 1$, $\alpha_f = 0$,
$M_1 = M_2 = M_f = 10^{8}$.
The calculation is performed by $T_{max} = 100$ with time step $\tau = 10$ for two - dimensional problem. For three - dimensional problem, we set $T_{max} = 6000$ with time step $\tau = 300$.
Initial condition is $p^0 = 0.1*10^7$ and $g = 0.2*10^7$ for pressure boundary condition. 

To compare the results, we use relative errors $L_2$ between multiscale solution and fine-scale solution 
\[
e^{p_i}_{L^2} =   \left(
\frac{\int_{\Omega} (p_i - p^{ms}_i)^2 dx}{\int_{\Omega} p^2_i dx} \right)^{1/2}, \quad  
e^u_{L^2} =  \left( 
\frac{\int_{\Omega} (u - u^{ms})^2 \, dx}{\int_{\Omega} u^2 \, dx} 
\right)^{1/2}, 
\]
where $i$ is the index of the continuum ($i=1,2$), $y^{ms} = (p^{ms}_1, p^{ms}_2, u^{ms})$ is the multiscale solution using GMsFEM and $y = (p_1, p_2, u)$ is the fine grid solution. 

We use $DOF_h$ (Degree of Freedom) to denote fine grid system size and $DOF_H$ to denote problem size of the coarse scale system using GMsFEM. On GMsFEM, we choose a same number of multiscale basis functions in each local domains $M = M_p = M_u$, where  $M_p$ and $M_u$ are the number of basis functions.  
We use GMSH software to construct computational domains and grids \cite{geuzaine2009gmsh}. The implementation is based on the open-source library FEniCS \cite{fenics}.

\subsection{Two-dimensional problem}

\begin{figure}[h!]
\centering
\includegraphics[width=0.3\linewidth]{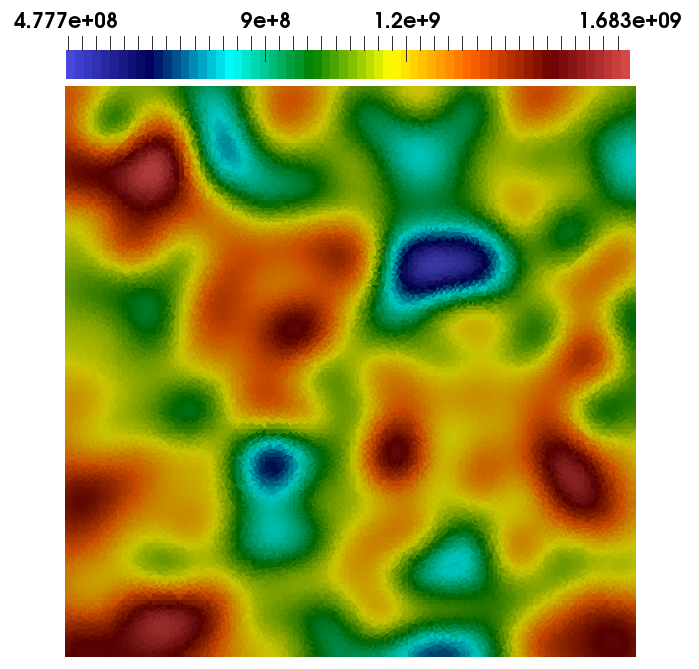}
\includegraphics[width=0.3\linewidth]{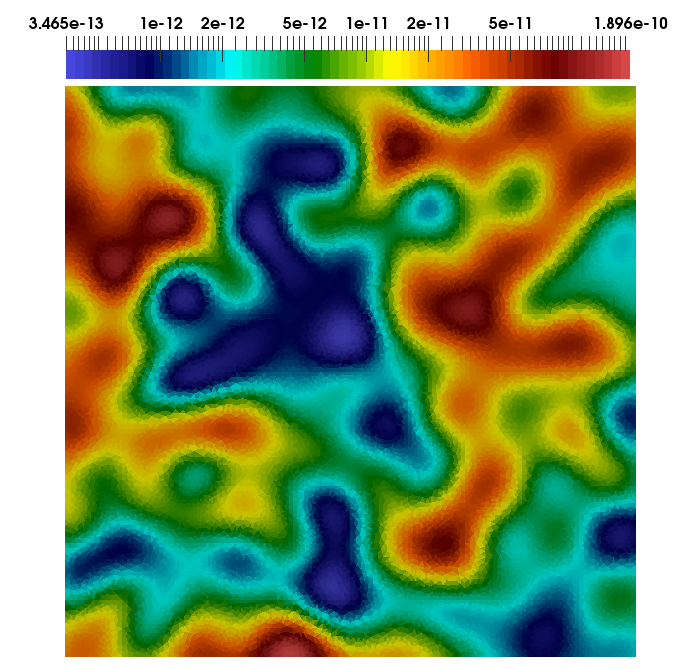}
\includegraphics[width=0.3\linewidth]{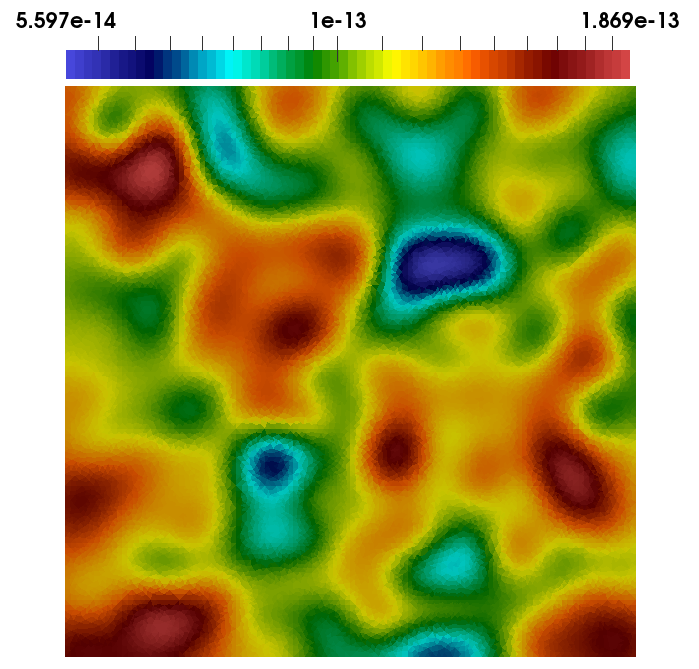}
\caption{Elasticity parameter $E$ (left) and heterogeneous permeabilities $k_1$(center) and $k_2$(right) for two - dimensional problem}
\label{kx2d}
\end{figure}

We simulate a two-dimensional model problem in heterogeneous fractured porous media. We set $k_f=10^{-4}$ and $r_{12} = 1000 \cdot k_2$.
Heterogeneous coefficients for elasticity modulus and heterogeneous permeability for first and second continuum are presented in Figure \ref{kx2d}.

\begin{figure}[h!]
\centering
\includegraphics[width=0.245\linewidth]{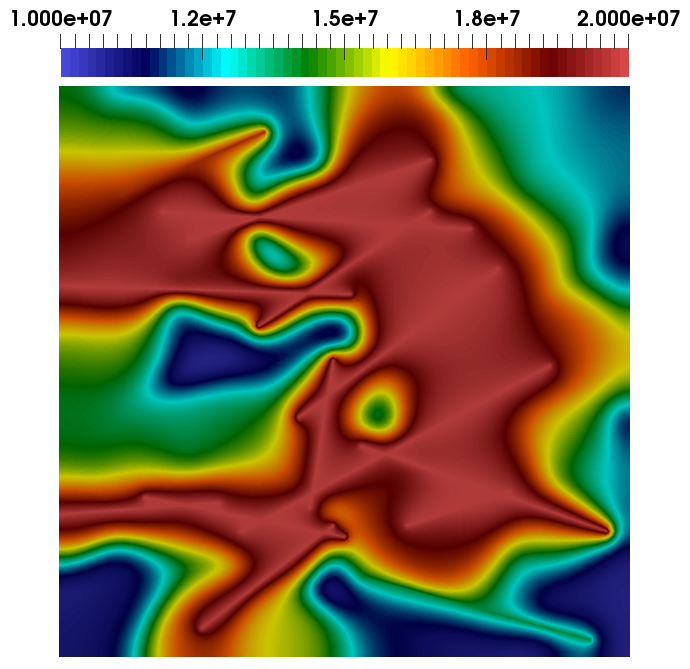}
\includegraphics[width=0.245\linewidth]{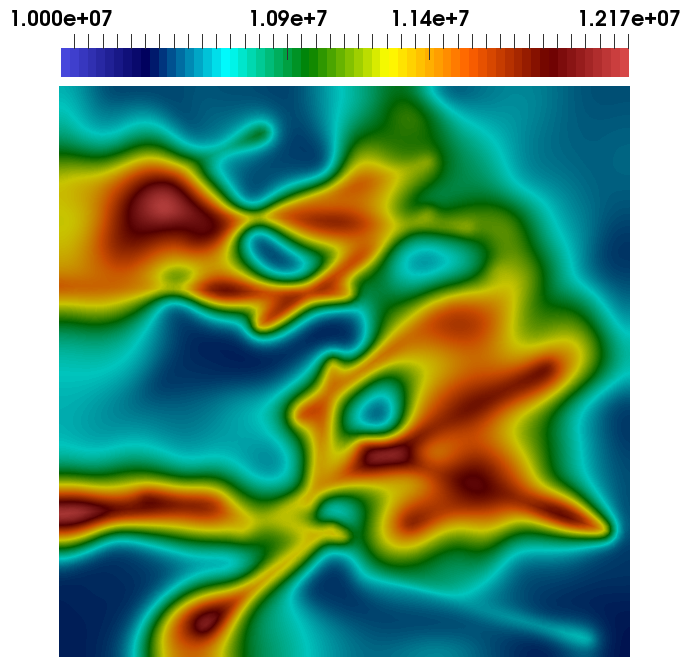}
\includegraphics[width=0.245\linewidth]{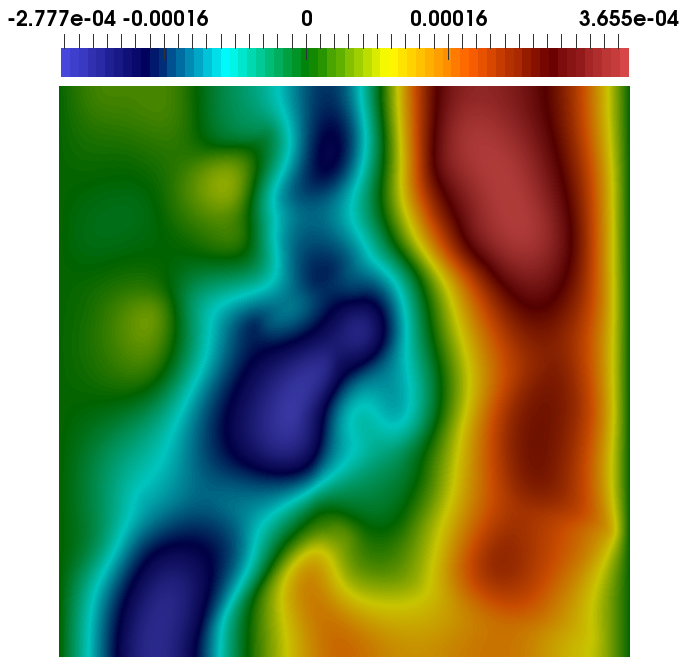}
\includegraphics[width=0.245\linewidth]{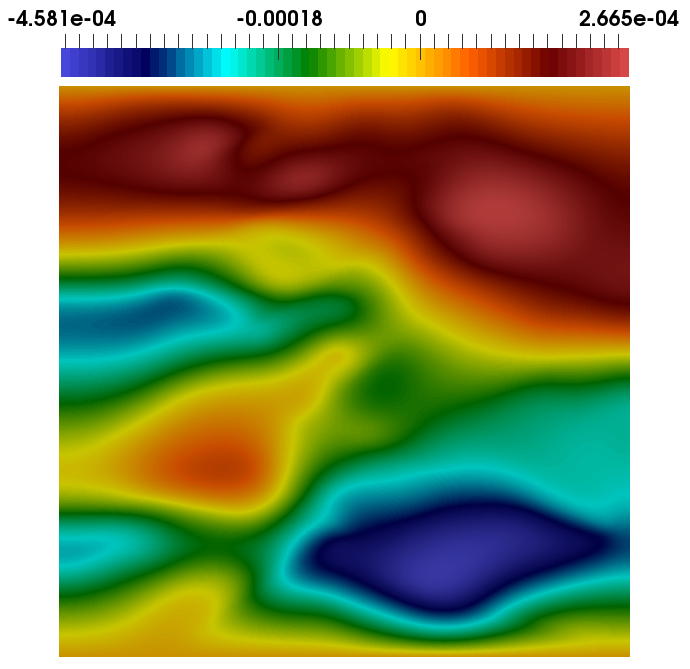}\\
\includegraphics[width=0.245\linewidth]{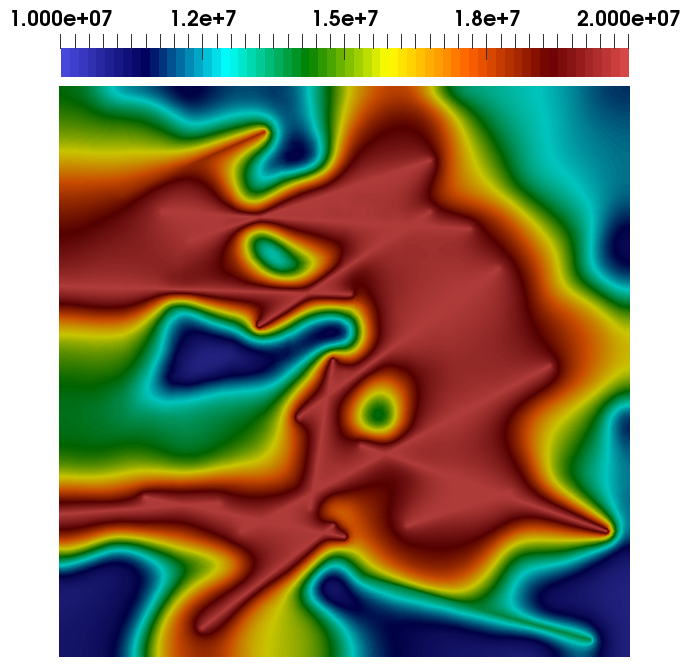}
\includegraphics[width=0.245\linewidth]{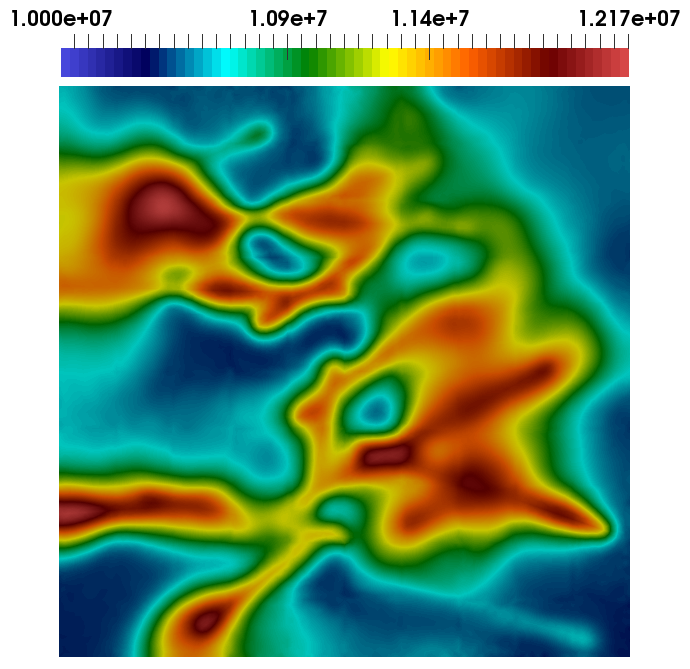}
\includegraphics[width=0.245\linewidth]{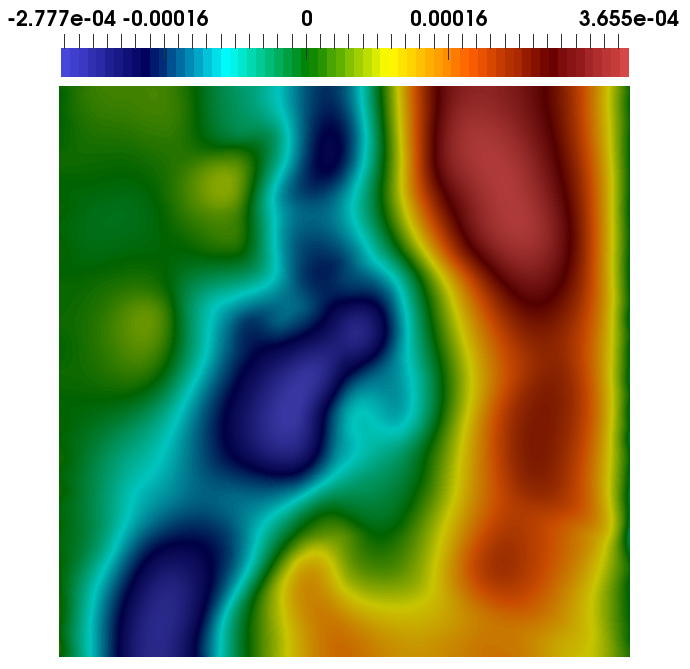}
\includegraphics[width=0.245\linewidth]{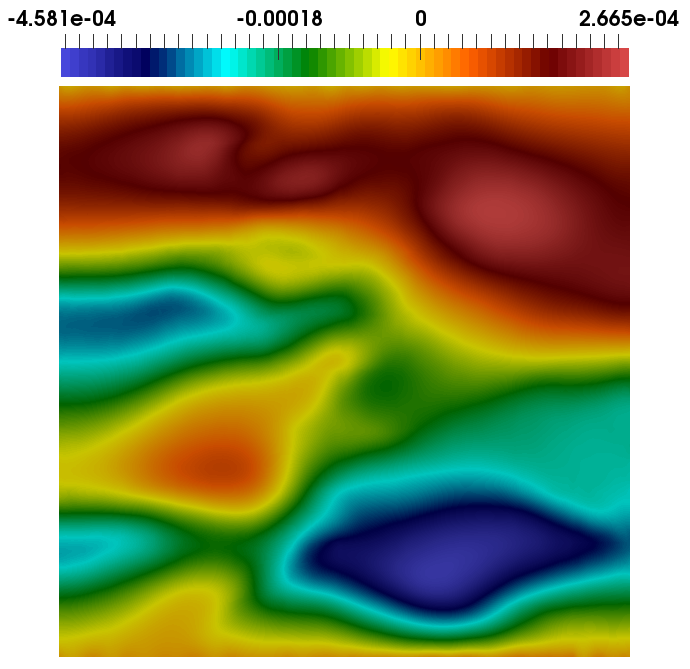}
\caption{Numerical results for two - dimensional problem. 
Pressure for first continuum and second continuum, displacement $X$ and $Y$ directions at final time (from left to right).
First row: fine scale solution. Second row: multiscale solution using 16 multiscale basis functions}
\label{u2d}
\end{figure}

\begin{table}[h!]
\begin{center}
\begin{tabular}{ | c | c | c | c | c | }
\hline
 \multirow{2}{*}{$M$}
 &  \multirow{2}{*}{$DOF_H$ }
 & \multicolumn{1}{c|}{$p_1$}
 & \multicolumn{1}{|c|}{$p_2$}
 & \multicolumn{1}{|c|}{$u$} \\
 & 
& $e^{p_1}_{L^2}$  (\%)  
& $e^{p_2}_{L^2}$  (\%)  
& $e^{u}_{L^2}$ (\%)  \\ 
\hline
\multicolumn{5}{|c|}{Coarse mesh $5 \times 5$} \\
\hline
1 		& 144 	& 28.501   & 9.514 	& 93.582   \\
2 		& 288 	&  24.552 	& 7.933 	& 90.760  \\
4 		& 576 	 & 10.228   & 3.258 	& 53.017 \\
8 		& 	1152	& 5.695   & 2.080 	& 27.283  \\
12 	& 1728	& 2.526 & 0.887 		& 12.646   \\
16 	& 2304	& 1.755  & 0.582 		& 9.076   \\
\hline
\multicolumn{5}{|c|}{Coarse mesh  $10 \times 10$} \\
\hline
1 		& 484 		 & 24.385   & 7.825 	& 111.247  \\
2 		& 968 		& 11.947   & 3.906  & 40.595   \\
4 		& 1936 	& 5.046   & 1.766 	& 21.797   \\
8 		& 	3872	& 1.172  & 0.384 	& 4.406  \\
12 	& 5808	& 0.270  & 0.149 		& 3.166   \\
16 	& 7744	& 0.154  & 0.083 		& 2.410  \\
\hline
\end{tabular}
\end{center}
\caption{Relative errors (\%) for displacement and pressures with different number of multiscale basis functions  for two - dimensional problem. Coarse grids $5 \times 5$  and $10 \times 10$. $DOF_h = 57504$}
\label{table2d}
\end{table}

\begin{table}[h!]
\begin{center}
\begin{tabular}{ | c | c | c | c | c | }
\hline
 \multirow{2}{*}{$M_u$}
 &  \multirow{2}{*}{$DOF_H$ }
 & \multicolumn{1}{|c|}{$p_1$}
 & \multicolumn{1}{|c|}{$p_2$} 
 & \multicolumn{1}{|c|}{$u$}\\
 & 
& $e^{p_1}_{L^2}$  (\%) 
& $e^{p_2}_{L^2}$  (\%) 
& $e^{u}_{L^2}$ (\%) \\ \hline
\multicolumn{5}{|c|}{\textbf{$M_p = 1$ }} \\
\hline
1 		& 484 	& 24.385 & 7.825  & 111.247   \\
\hline
\multicolumn{5}{|c|}{\textbf{$M_p = 2$ }} \\
\hline
1 		& 726 	& 11.905  & 3.741 	& 72.299  \\
2 		& 968 	& 11.947  & 3.906 	& 40.595  \\
\hline
\multicolumn{5}{|c|}{\textbf{$M_p = 4$ }} \\
\hline
1 		& 1210 	& 5.029   & 1.723 & 51.183   \\
2 		& 1452 	& 5.046  & 1.772 & 23.402  \\
4 		& 1936 	& 5.046   & 1.766 & 21.797  \\
\hline
\multicolumn{5}{|c|}{\textbf{$M_p = 8$ }} \\
\hline
1 		& 2178 	& 1.169   & 0.585 & 38.968   \\
2 		& 2420 	& 1.170   & 0.426 & 9.469   \\
4 		& 2904 	& 1.170   & 0.406 & 8.285   \\
8 		& 3872	& 1.172  & 0.384 & 4.406  \\
\hline
\multicolumn{5}{|c|}{\textbf{$M_p = 12$ }} \\
\hline
1 		& 3146 	& 0.302  & 0.504 	& 38.068   \\
2 		& 3388 	& 0.272  & 0.237 	& 8.940   \\
4 		& 3872 	& 0.271  & 0.194 	& 7.705   \\
8 		& 4840	& 0.271  & 0.154 	& 3.814  \\
12 	        & 5808	& 0.270  & 0.149 	& 3.166  \\
\hline
\multicolumn{5}{|c|}{\textbf{$M_p = 16$ }} \\
\hline
1 		& 4114 	& 0.212  & 0.485 	& 37.974   \\
2 		& 4356 	& 0.156  & 0.196 	& 9.044   \\
4 		& 4840 	& 0.155 & 0.148 	& 7.775   \\
8 		& 5808	& 0.155 & 0.095 	& 3.836  \\
12 	        & 6776	& 0.155  & 0.087 	& 3.211   \\
16 	        & 7744	 & 0.154  & 0.083 	& 2.410  \\
\hline
\end{tabular}
\end{center}
\caption{Relative errors (\%) for displacement and pressures with different number of multiscale basis functions  for two - dimensional problem. Coarse grid $10 \times 10$. $DOF_h = 57504$}
\label{table2d2}
\end{table}

In Figures \ref{u2d} distribution of pressure for first continuum and second continuum, displacement along $X$ and $Y$ directions at final time are presented. 
On the first row, we depict a fine scale solution and multiscale solution with 16 multiscale basis functions for GMsFEM is presented on second row. 
Comparing the fine-scale solution with the multiscale solution with 16 basis functions in Figure \ref{u2d} for displacement along $X$ and $Y$ pressure, we can observe good accuracy. In general, solutions look good without visible oscillations. 

  In Table  \ref{table2d}, we present an errors for $5 \times 5$  and $10 \times 10$ coarse grid with $M = M_p = M_u$. 
The results  show that 8 multiscale basis functions are enough to achieve good results with $4.406 \%$ of $L^2$ error for displacement, $1.172 \%$ and  $0.384 \%$ of  $L_2$ errors for first continuum and second continuum pressures. 
When we increase number of the multiscale basis functions, the relative $L^2$ errors are decrease two times when we take two multiscale basis functions instead of one. We have similar improvements for further increasing of the multiscale basis functions.  

Table \ref{table2d2} shows a comparison of the difference between multiscale and fine grid solutions, where we present relative $L_2$  for different number of the multiscale basis functions on the coarse grid $10 \times 10$. 
From the Table \ref{table2d2}, we observe that number of basis functions for pressure $M_p$ highly impact to the displacements errors. For example in the case $M_u = 2$, we have $40 \%$ of $L^2$ displacements error when we have 2 multiscale basis functions for pressure, and reduce error to $9 \%$ for $M_p = 8$. 
In conclusion, according to the results of the comparison, we observe a good convergence of the method when we increase number of the multiscale basis functions.

\subsection{Three-dimensional problem}

\begin{figure}[h!]
\centering
\includegraphics[width=0.31\linewidth]{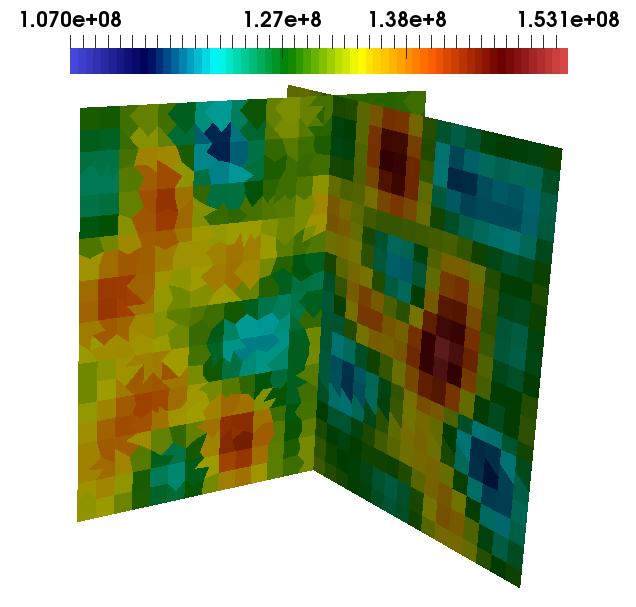} 
\includegraphics[width=0.31\linewidth]{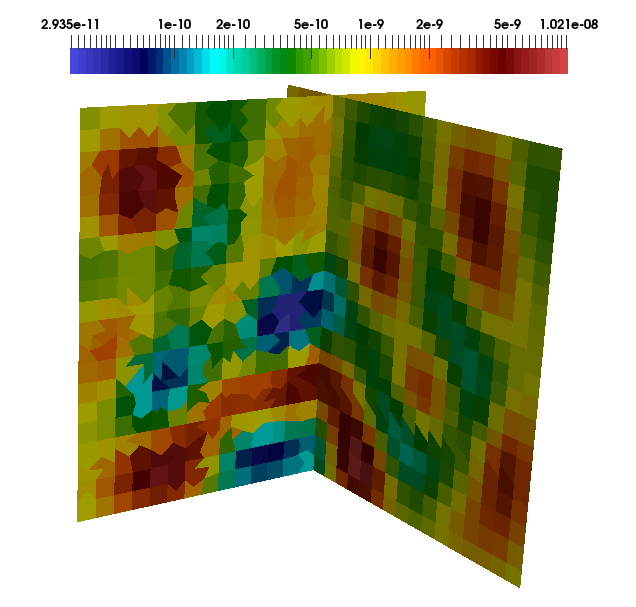} 
\includegraphics[width=0.31\linewidth]{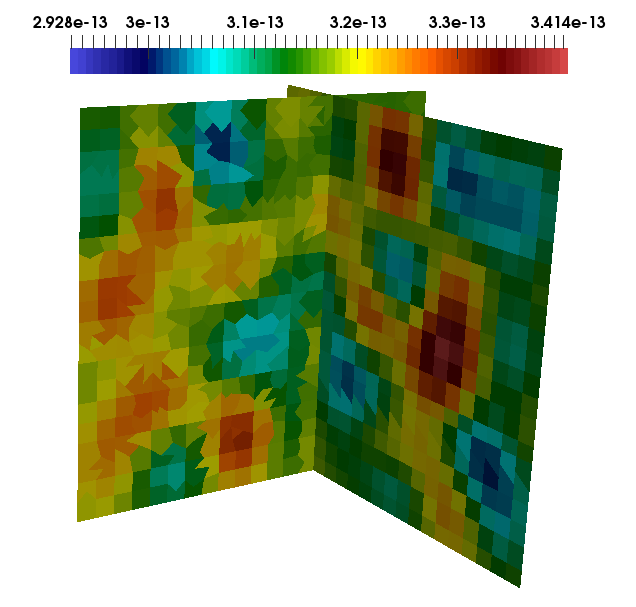}
\caption{Elasticity parameter $E$ (left) and heterogeneous permeabilities $k_1$(center) and $k_2$(right) for three - dimensional problem}
\label{kx3d}
\end{figure}

We simulate a three-dimensional model problem in heterogeneous fractured porous media. We set $k_f=10^{-4}$ and $r_{12} = 25 \cdot k_2$.
Heterogeneous coefficients for elasticity modulus and heterogeneous permeability for first continuum and second continuum are presented in Figure \ref{kx3d}.

\begin{figure}[h!]
\centering
\vspace{2.0 mm}
\includegraphics[width=0.31\linewidth]{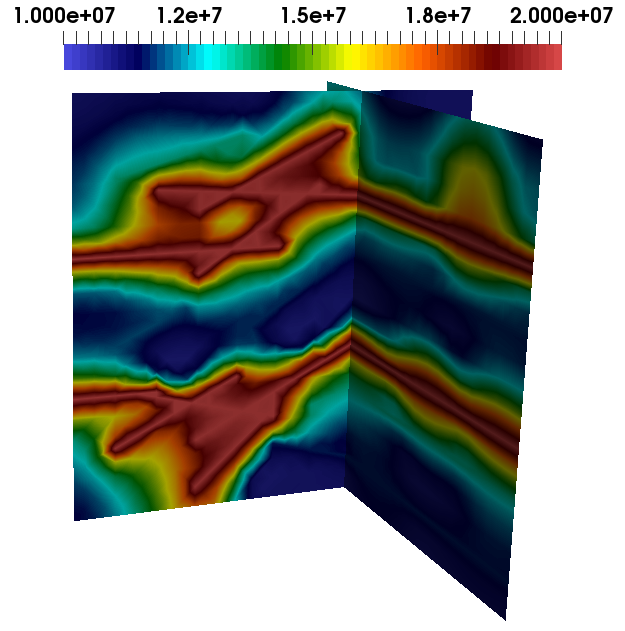}
\includegraphics[width=0.31\linewidth]{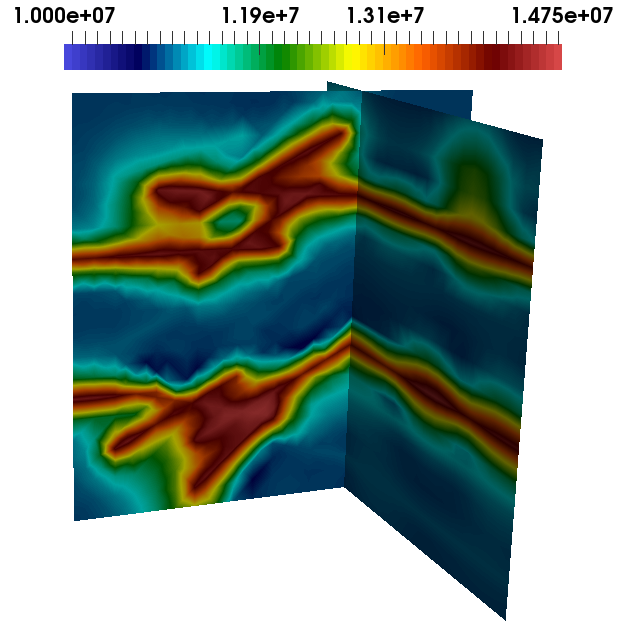}  \\
\includegraphics[width=0.31\linewidth]{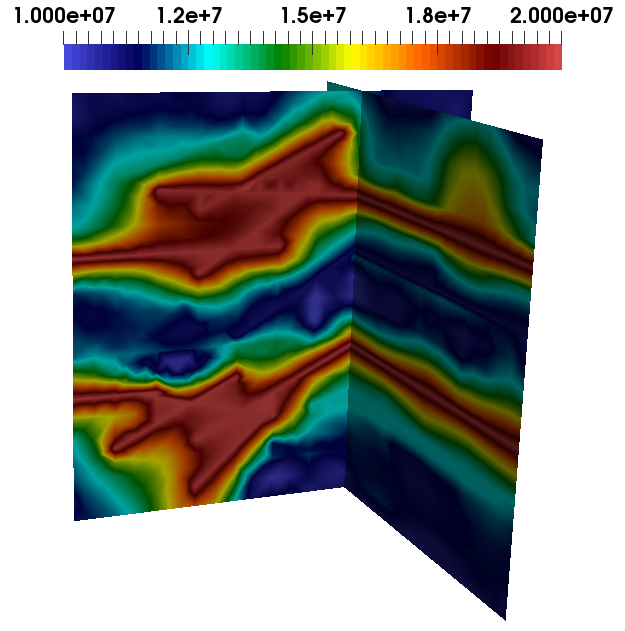}
\includegraphics[width=0.31\linewidth]{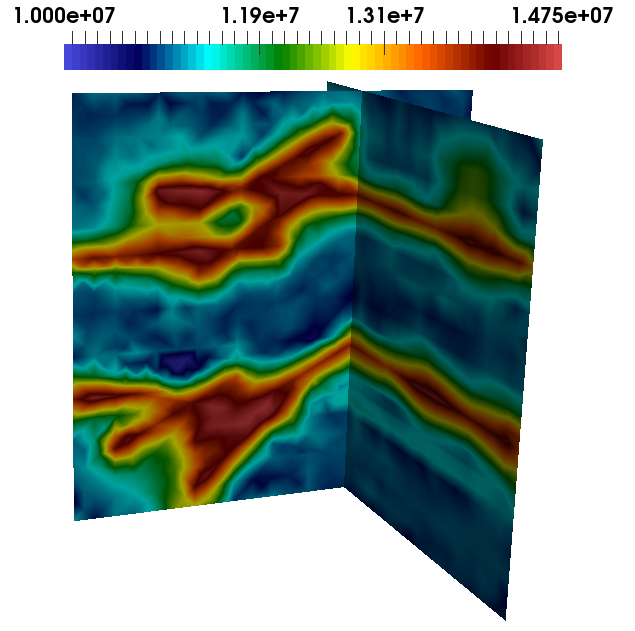} 
\caption{Numerical results for three - dimensional problem. 
Pressure for first continuum and second continuum (from left to right).
First row: fine scale solution. Second row: multiscale solution using 16 multiscale basis functions}
\label{p3d}
\end{figure}

\begin{figure}[h!]
\centering
\vspace{2.0 mm}
\includegraphics[width=0.31\linewidth]{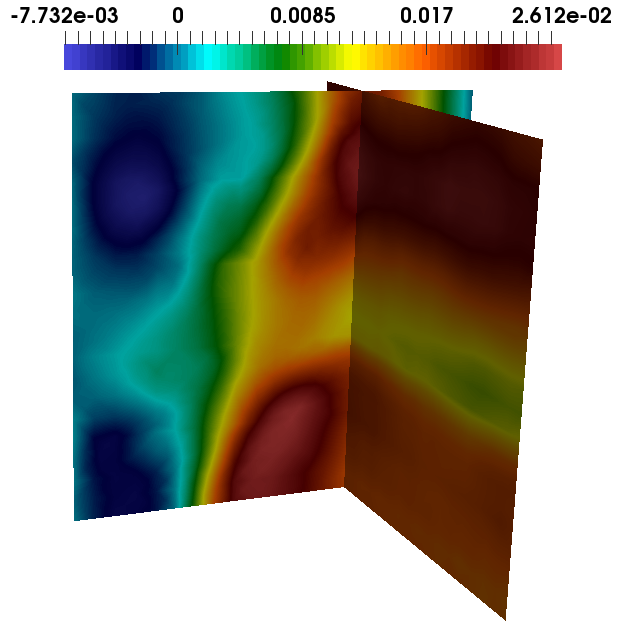} 
\includegraphics[width=0.31\linewidth]{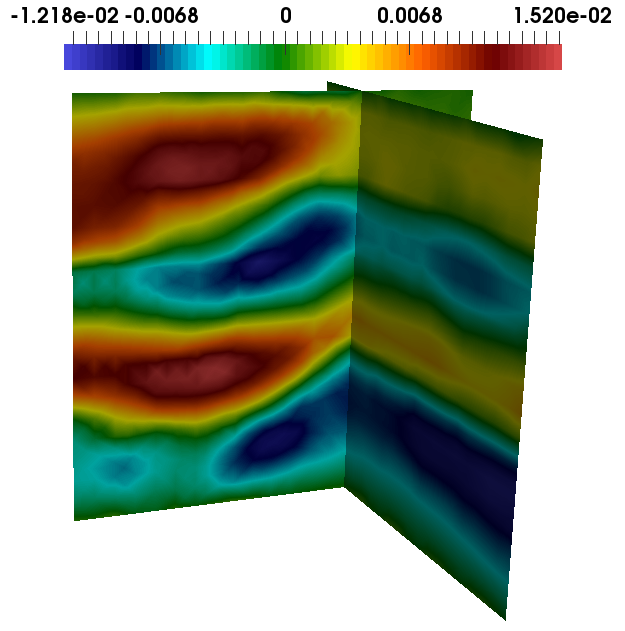}
\includegraphics[width=0.31\linewidth]{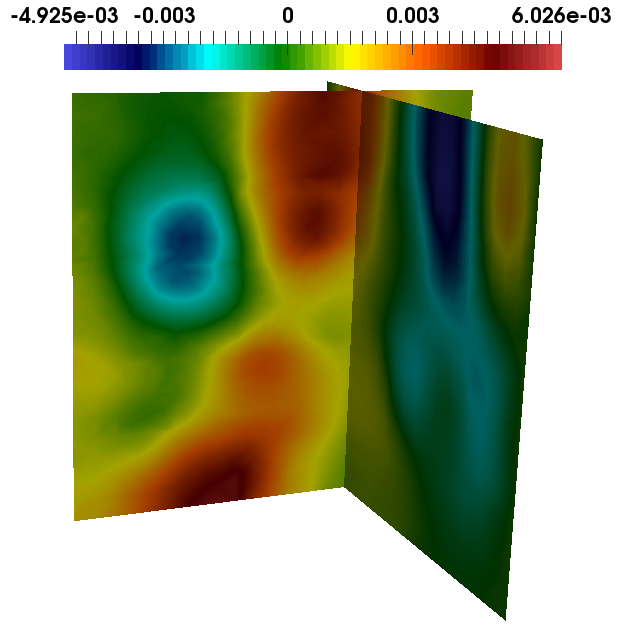}
\includegraphics[width=0.31\linewidth]{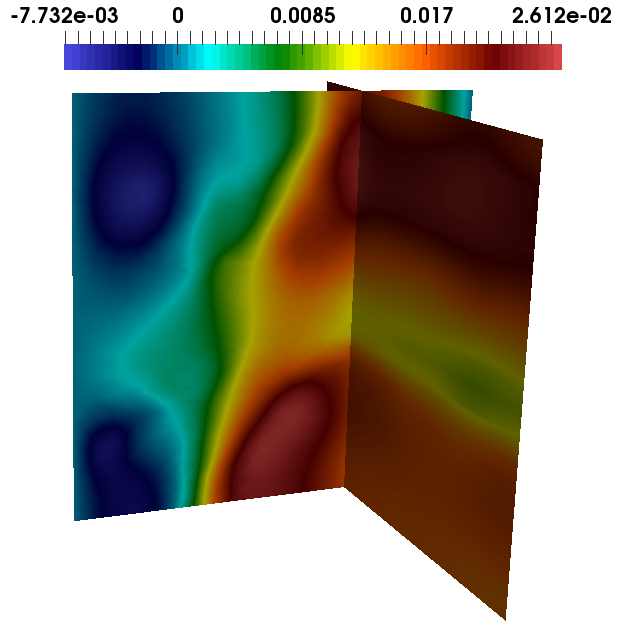} 
\includegraphics[width=0.31\linewidth]{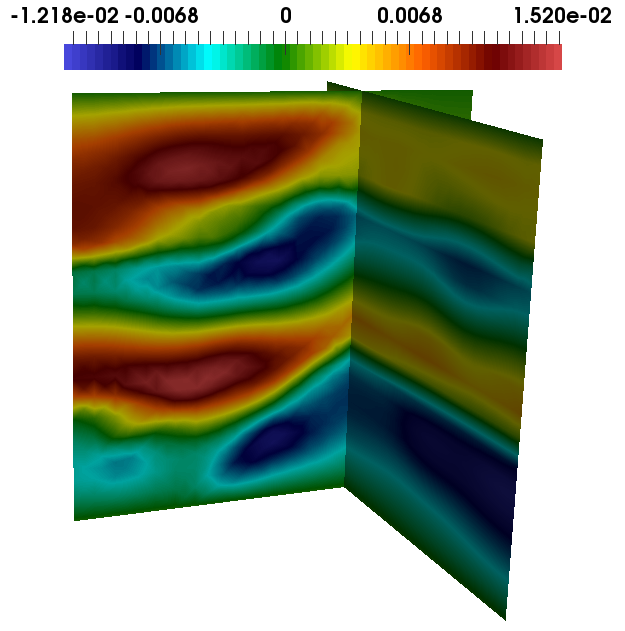}
\includegraphics[width=0.31\linewidth]{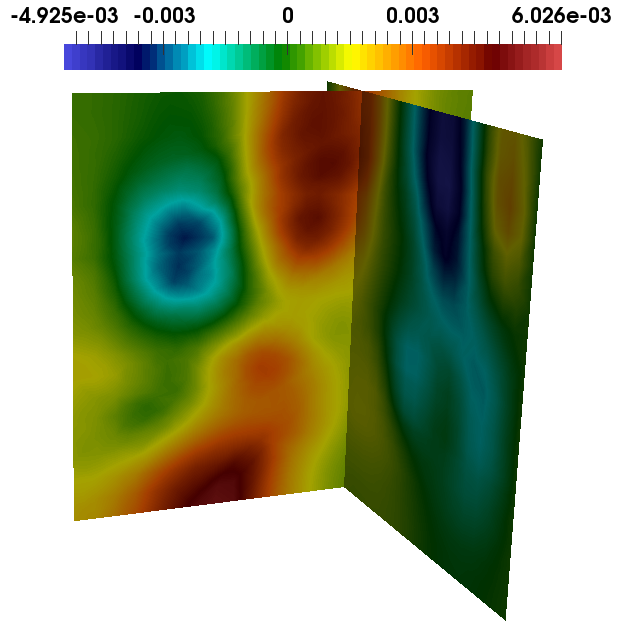}
\caption{Numerical results for three - dimensional problem. 
Displacement $X$, $Y$ and $Z$ directions at final time (from left to right).
First row: fine scale solution. Second row: multiscale solution using 16 multiscale basis functions}
\label{u3d}
\end{figure}

\begin{table}[h!]
\begin{center}
\begin{tabular}{ | c | c | c | c | c | }
\hline
 \multirow{2}{*}{$M_u$}
 &  \multirow{2}{*}{$DOF^u_H$ }
 & \multicolumn{1}{|c|}{$p_1$}
 & \multicolumn{1}{|c|}{$p_2$}
 & \multicolumn{1}{|c|}{$u$} \\
 & 
& $e^{p_1}_{L^2}$  (\%)   
& $e^{p_2}_{L^2}$  (\%)  
& $e^{u}_{L^2}$ (\%)    \\ \hline
\multicolumn{5}{|c|}{\textbf{$M_p = 1$ }} \\
\hline
1 	& 1080 & 25.713   & 29.035  & 92.385  \\
\hline
\multicolumn{5}{|c|}{\textbf{$M_p = 2$ }} \\
\hline
1 		& 1512  & 20.461   & 22.644 	& 130.73   \\
2 		& 2106 & 20.496   & 22.648 	& 82.179   \\
\hline
\multicolumn{5}{|c|}{\textbf{$M_p = 4$ }} \\
\hline
1 		& 2376 	& 6.453   & 6.991 	& 86.899   \\
2 		& 3024 	& 6.557   & 7.114 	& 59.549   \\
4 		& 4320 	& 6.626   & 7.190 	& 22.580   \\
\hline
\multicolumn{5}{|c|}{\textbf{$M_p = 8$ }} \\
\hline
1 		& 4104 	& 3.563   & 3.869 	& 78.466   \\
2 		& 4752 	& 3.647   & 3.974 	& 44.919   \\
4 		& 6048 	& 3.711   & 4.050 	& 12.562  \\
8 		& 8640  & 3.709   & 4.052 	& 11.773  \\
\hline
\multicolumn{5}{|c|}{\textbf{$M_p = 12$ }} \\
\hline
1 		& 5832 	& 2.444   & 2.710 	& 77.902   \\
2 		& 6480 	& 2.481   & 2.755 	& 43.286   \\
4 		& 7776 	& 2.530   & 2.813 	& 10.389   \\
8 		& 10368	& 2.529   & 2.813 	& 8.518  \\
12 	        & 12960	& 2.530   & 2.814 	& 8.518   \\
\hline
\multicolumn{5}{|c|}{\textbf{$M_p = 16$ }} \\
\hline
1 	& 7560 	& 1.898   & 2.080 	& 76.300   \\
2 	& 8208 	& 1.897   & 2.077 	& 41.473   \\
4 	& 9504 	& 1.889   & 2.064 	& 9.687   \\
8 	& 12096	& 1.888   & 2.063 	& 6.908  \\
12 	& 14688	& 1.888   & 2.064 	& 6.170   \\
16 	& 17280	& 1.887   & 2.063 	& 5.881  \\
\hline
\end{tabular}
\end{center}
\caption{Relative errors for displacement and pressures with different number of multiscale basis functions  for three - dimensional problem. $DOF_h = 108045$}
\label{table3d2}
\end{table}

In Figures \ref{p3d}-\ref{u3d} distribution of pressure for first continuum and second continuum, displacement along $X$,$Y$ and $Z$ directions at the final time are presented. On the first row, we depict a fine scale solution and multiscale solution with 16 multiscale basis functions for GMsFEM is presented on second row. We observe a good results of the multiscale method compared with the fine grid solution.

Table \ref{table3d2} shows a comparison of the difference between multiscale and fine grid solutions. We present relative $L^2$ for displacements and pressures. For three - dimensional model problem using 8 multiscale basis functions, we have
$11.773 \%$ of $L^2$ error for displacement, $3.709 \%$ and $4.052 \%$ of $L_2$ errors for first continuum and second continuum pressures. When we take 16 multiscale basis functions, we obtain two times better results with $5.881\%$ of $L^2$ error for displacement, $1.887 \%$ and $2.063 \%$ of $L^2$ errors for first continuum  and second continuum pressures.
For the three-dimensional problem, we also observe good convergences for the poroelasticity problem in heterogeneous and fractured media.

\section{Conclusion}

In this work, we considered the poroelasticity problem in heterogeneous and fractured medium. We presented a mathematical model and fine grid approximation using finite element method for general multicontinuum poroelasticity problem. We developed a generalized multiscale finite element framework and explained construction of the multiscale basis functions. 
We presented results of the numerical investigation for model problems in two and three-dimensional formulations. We compared a relative error between multiscale and fine-scale solutions for different number of multiscale basis functions. The proposed multiscale method provides a good accuracy for two and three-dimensional problems in heterogeneous fractures media
with few degrees of freedoms.

\section{Acknowledgements}
AT's, MV's and DS's works are supported by the mega-grant of the Russian Federation Government N14.Y26.31.0013 and RSF N17-71-20055.
The research of Eric Chung is partially supported by the Hong Kong RGC General Research Fund (Project numbers 14304217 and 14302018) and CUHK Faculty of Science Direct Grant 2018-19.

\bibliographystyle{unsrt}
\bibliography{lit}

\end{document}